\documentclass[12pt]{amsart}
\usepackage{latexsym, amssymb, amsmath, amscd}
 \usepackage{eucal}

    \newtheorem{rema}{Remark}[section]

   \newtheorem{theo}[rema]{Theorem}
   \newtheorem{def-theo}[rema]{Definition-Theorem}
 \newtheorem{conj}[rema]{Conjecture}
   
    \newtheorem{lemma}[rema]{Lemma}
    \newtheorem{corol}[rema]{Corollary}
     \newtheorem{exam}[rema]{Example}
  \newtheorem{rmk}[rema]{Remark}

	\newcommand{\nno}{\nonumber}

	\newcommand{\p}{\partial}

 \newcommand{\pf}{{\it Proof:}\hspace{2ex}}
 
 \newcommand{\epfv}{\hspace{1em}$\Box$\vspace{1em}}

\newcommand{\bC}{{\mathbb C}}
\newcommand{\bZ}{{\mathbb Z}}
\newcommand{\bR}{{\mathbb R}}
\newcommand{\bQ}{{\mathbb Q}}
\newcommand{\bN}{{\mathbb N}}

\newcommand{\cL}{{\mathcal L}}

\newcommand{\cA}{{\mathcal A}}

\newcommand{\mcD}{{\mathcal D}}

\newcommand{\bD}{{\mathbb D}}

\newcommand{\cM}{{\mathcal M}}

\newcommand{\cB}{{\mathcal B}}

\newcommand{\cC}{{\mathcal C}}

\newcommand{\cE}{{\mathcal E}}

\newcommand{\mcZ}{{\mathcal Z}}
\newcommand{\mcE}{{\mathcal E}}

\newcommand{\cz}{\bC[z]}
\newcommand{\czz}{\bC[z^{-1}, z]}

\newcommand{\cxz}{\bC[\xi, z]}
\newcommand{\cxi}{\bC[\xi]}
\newcommand{\Kz}{K[z]}
\newcommand{\cAz}{\cA[z]}

\newcommand{\Ker}{\rm{Ker\,}}
\newcommand{\im}{\operatorname{Im}}

\newcommand{\JC}{{\bf JC }}
\newcommand{\IC}{{\bf IC }}
\newcommand{\VC}{{\bf VC }}

\newcommand{\BQ}{\begin{eqnarray}}
\newcommand{\EQ}{\end{eqnarray}}
\newcommand{\BQn}{\begin{eqnarray*}}
\newcommand{\EQn}{\end{eqnarray*}}

\title[Images of Commuting Differential Operators]
{Images of Commuting  Differential Operators of Order One with Constant Leading Coefficients}

  \author{Wenhua Zhao}      

\begin{document}

\begin{abstract}
We first study some properties of images of commuting differential operators of polynomial algebras of order one with constant  leading coefficients. We then propose what we call the {\it image} conjecture on these differential operators and show that the 
{\it Jacobian} conjecture \cite{BCW}, \cite{E}, 
\cite{Bo} (hence also the {\it Dixmier} conjecture \cite{D}) and the {\it vanishing}  conjecture \cite{GVC} of differential operators with constant coefficients are actually equivalent to certain special cases of the {\it image} conjecture. A connection of the {\it image} conjecture, and hence also the {\it Jacobian} conjecture, with multidimensional Laplace transformations of polynomials is also discussed. 
\end{abstract}




\thanks{The author has been partially supported 
by NSA Grant R1-07-0053}

 \bibliographystyle{alpha}
    \maketitle


\renewcommand{\theequation}{\thesection.\arabic{equation}}
\renewcommand{\therema}{\thesection.\arabic{rema}}
\setcounter{equation}{0}
\setcounter{rema}{0}
\setcounter{section}{0}

\section{\bf Introduction}

\subsection{Background and Motivation}
Let $z=(z_1, z_2, \dots, z_n)$ 
be $n$ commutative free variables and $\cz$ 
the polynomial algebra in $z$ over $\bC$. 
Recall that the Jacobian conjecture ({\bf JC})
proposed by O. H. Keller \cite{Ke} in 1939  
claims that {\it any polynomial map $F$ 
of $\bC^n$ with Jacobian $j(F)\equiv 1$ must 
be an automorphism of $\bC^n$}. 
Despite intense study 
from mathematicians in the last seventy years, 
the conjecture is still open
even for the case $n=2$. In 1998, the Fields Medalist 
and also the Wolf Prize Winner S. Smale \cite{Sm} 
included the Jacobian conjecture 
in his list of $18$ fundamental mathematical problems 
for the $21$st century. For more history and 
known results on the Jacobian conjecture, 
see \cite{BCW}, \cite{E}, \cite{Bo} and 
references therein.

Recently, the equivalence of the \JC with   
the Dixmier conjecture proposed by J. Dixmier \cite{D} in 
$1968$ has been established first by Y. Tsuchimoto 
\cite{Ts} in $2005$ and later by A. Belov and M. Kontsevich
\cite{BK} and P. K. Adjamagbo and 
A. van den Essen \cite{AE} in $2007$. 
The implication of the Jacobian conjecture 
from the Dixmier conjecture was actually proved 
much earlier by V. Kac (unpublished 
but see \cite{BCW}) in $1982$. 

The Dixmier conjecture claims that {\it any homomorphism of the Weyl algebras must be an automorphism of the Weyl algebras}.
Note that in \cite{AE} P. K. Adjamagbo and A. van den Essen 
also showed that the {\bf JC} is also equivalent to what they called the {\it Poisson} conjecture on Poisson algebras.

This paper is mainly motivated by the connections discussed  
in \cite{HNP} and \cite{GVC} of the \JC with the  
{\it vanishing} conjecture ({\bf VC})  on  
differential operators (of any order) with constant 
coefficients of the polynomial algebras $\cz$.

First let us recall  the {\it vanishing} conjecture 
({\bf VC}) proposed in \cite{HNP} and \cite{GVC}. 
For later purposes, here we put it in 
a more general form. 

\begin{conj} \label{VC} $(${\bf The Vanishing Conjecture}$)$
Let $P(z), Q(z)\in \cz$ and $\Lambda$ be a differential operator 
 of $\cz$ with constant coefficients. Assume that  
$\Lambda^m(P^m)=0$ for each $m\ge 1$. Then we have 
$\Lambda^m(P^m Q)=0$ when $m\gg 0$.
\end{conj}

The motivation for the conjecture above is the 
following theorem proved in \cite{HNP}.

\begin{theo}\label{VC2JC}
Let $\Delta=\sum_{i=1}^n \p^2/\p z_i^2$, i.e.\@ the Laplace operator of $\cz$. Then 
the \JC holds for all $n\ge 1$ iff 
the \VC holds for all $n\ge 1$ with $\Lambda=\Delta$ 
and $P(z)=Q(z)$, where $P(z)$ is any homogeneous 
polynomial in $z$ of degree $4$.
\end{theo}

The proof of the theorem above is based on the classical celebrated homogeneous reduction of the \JC achieved by H. Bass, E. Connell, D. Wright \cite{BCW} and A. V. Yagzhev \cite{Y} and also based on the remarkable symmetric reduction on the \JC achieved independently by M. de Bondt and A. van den Essen \cite{BE} and G. Meng \cite{Me}. For more details of the proof, see \cite{BurgersEq} and \cite{HNP}. 

Note that it has been shown 
in \cite{EZ} that Theorem \ref{VC2JC} 
also holds without the condition $P(z)=Q(z)$ and in \cite{GVC} 
that it also holds if the Laplace operators $\Delta$ are replaced by any sequence of quadratic homogeneous 
differential operators 
$\{\Lambda_n (\p)\,|\, n\in \bN\}$ with the rank of 
the quadratic form $\Lambda_n(\xi)$ goes to $\infty$ 
as $n\to \infty$, where 
$\Lambda_n(\xi)$ is the symbol of 
the differential operator  
$\Lambda_n(\p)$ $(n\in \bN)$.

In this paper, we will show that the \JC (hence also the 
Dixmier conjecture \cite{D} and the Poisson conjecture 
in \cite{AE}) and the \VC above are actually equivalent 
to certain special cases of what we call the {\it image} 
conjecture (See Conjecture \ref{IC} below) of commuting differential operators of $\cz$ of order one with constant leading coefficients. See Theorems \ref{VC2IC}, \ref{JC2IC} and \ref{AJC2IC} for the precise statements. To formulate the image conjecture ({\bf IC}), let us first fix some notations that will be used throughout this paper.

Let $\cA$ be any commutative algebra over a field of characteristic zero and $\cAz$ the 
polynomial algebra of $z$ over $\cA$. 
For any $1\le i\le n$, we set $\p_i\!:=\p/\p z_i$ and 
$\p\!:=(\p_1, \p_2, ..., \p_n)$. 
We denote by $\mcD_\cA[z]$ the Weyl algebra 
of differential operators of $\cAz$, and 
$\bD_\cA[z]$ the subspace of differential operators of $\cAz$ 
of order one with constant leading coefficients, i.e.\@  the differential operators of the form 
$h(z)+\sum_{i=1}^n c_i\p_i$ for some $h(z)\in \cAz$
and $c_i\in \cA$ $(1\le i\le n)$. 

For any $\cC \subset \mcD_\cA[z]$, we define the 
{\it image}, denoted by $\im\cC$, of $\cC$ to be 
$\sum_{\Phi\in \cC} \Phi (\cAz)$. Furthermore, 
we say $\cC$ is {\it commuting} if any two differential 
operators in $\cC$ commute with each other. 

\begin{conj} \label{IC} $(${\bf The Image Conjecture}$)$
Let $\cA$ be any commutative algebra over a field of characteristic zero and $\cC$ a commuting subset of differential operators of $\cAz$ of order one with constant leading coefficients. Then for any $f, g\in \cAz$ with $f^m\in \im\cC$ for each $m\ge 1$, 
we have $f^m g\in \im\cC$ when $m\gg 0$.
\end{conj}

Note that the statement of the {\it image} conjecture 
({\bf IC}) above is similar to the statement of the Mathieu conjecture proposed by O. Mathieu \cite{Ma} in $1995$. 

\begin{conj} \label{MC} $(${\bf The Mathieu Conjecture}$)$
Let $G$ be a compact connected real Lie group with the 
Haar measure $d\sigma$. Let $f$ be a complex-valued $G$-finite 
function over $G$ such that $\int_G f^m \, d\sigma=0$ 
for each $m\ge 1$. Then for any $G$-finite 
function $g$ over $G$, 
$\int_G f^m g \,  d \sigma =0$ when $m\gg 0$.
\end{conj}

In \cite{Ma} Mathieu also proved that his conjecture implies the 
{\bf JC}. Later J. Duistermaat and W. van der Kallen 
\cite{DK} proved the Mathieu conjecture 
for the case of tori, which can be re-stated as follows. 

\begin{theo}\label{ThmDK} 
$(${\bf Duistermaat and van der Kallen}$)$ 
Let $\czz$ be the algebra of Laurent polynomials over $\bC$,  
and $\cM \subset \czz$ the subspace of all Laurent polynomials 
in $z$ with no constant terms. Then for any 
$f(z), g(z)\in \czz$ with $f^m(z)\in \cM$ for 
each $m\ge 1$, we have $f^m(z) g(z)\in \cM$ 
when $m\gg 0$.
\end{theo}

Actually, if we consider certain generalizations of the \IC for differential operators of some localizations of $\cz$, especially, for those commuting differential operators of order one with leading coefficients related with classical orthogonal polynomials, the \IC is indeed connected with the Mathieu conjecture and the Duistermaat-van der Kallen theorem. For example, let $\Phi_i=\p_i-z_i^{-1}$ $(1\le i\le n)$ and $\cC=\{\Phi_1, \Phi_2, ..., \Phi_n\}$. Note that $\cC$ is a commuting subset of differential operators of the Laurent 
polynomial algebra $\czz$. Then one can show that the image  
$\im\cC=\sum_{i=1}^n\Phi_i(\czz)$ is exactly the subspace  
$\cM$ in Theorem \ref{ThmDK}. Hence, the \IC for $\czz$ and the differential operators $\cC$ above is equivalent to the Duistermaat-van der Kallen theorem. But, unfortunately, the straightforward generalization of the \IC to $\czz$ does not always hold. For a modified generalization of the \IC and their connections with the Mathieu conjecture and the Duistermaat-van der Kallen theorem, 
see \cite{GIC}.

\subsection{Arrangement} 
In section \ref{S2}, we assume that $\cA$ is a commutative algebra over a field $K$ of characteristic zero and prove some general properties of the images of commuting differential operators of order one with constant leading coefficients 
in $K$. We show in Theorem \ref{T2.1.4} that in this case, 
the \IC can be reduced to the case that  
$\cC=\{\Phi_1, \Phi_2, ..., \Phi_n\}$, where 
$\Phi_i=\p_i-\p_i(q(z))$ for some $q(z)\in \cAz$. 
We then assume that $\cA=K$ and 
$\cC$ as above except with $n$ replaced by any 
$1\le k\le n$, and show 
in Theorem \ref{T2.1.8} that with a properly defined 
$\mcD$-module structure on $K[z]$, $\im \cC$ forms 
a $\mcD$-submodule of $K[z]$, and the quotient 
$K[z]/ \im \cC$ is a holonomic $\mcD$-module. 
Consequently, we have that, when $k=n$, 
$\im\cC$ is a finite co-dimensional 
subspace of $K[z]$ 
(See Corollary \ref{C2.1.9}).

In Section \ref{S3}, we mainly consider the relations of 
the \VC and the \JC with the {\bf IC}. 
In Subsection \ref{S3.1}, we study 
the image of the commuting differential operators
$\Theta_i\!:=\xi_i-\p_i$ $(1\le i\le n)$ of $\cxz$, 
where $\xi=(\xi_1, \xi_2, ..., \xi_n)$ are another 
$n$ free commutative variables which also 
commute with $z$. 
Let $\Theta=\{\Theta_1, \Theta_2, ..., \Theta_n\}$. 
The main results of this subsection 
are Theorems \ref{Main-Thm-1} and \ref{D-Taylor} 
which identify $\im\Theta$ with the kernel of 
another linear map $\mcE:\cxz\to \cz$ 
(See Eq.\,(\ref{NewDef-E})). 
The kernel of $\mcE$ captures the condition 
that $\Lambda^m (f^m)=0$ $(m\ge 1)$ as well as 
the claim $\Lambda^m (f^mg)=0$ $(m\gg 0)$ in the {\bf VC}.
Therefore Theorems \ref{Main-Thm-1} 
and \ref{D-Taylor} provide a bridge connecting 
the {\bf VC}, hence also the {\bf JC}, the Dixmier Conjecture 
and the Poisson Conjecture, with the {\bf IC}. 

In Subsection \ref{S3.2}, by using the results obtained in 
Subsection \ref{S3.1}, we show that the \VC and the \JC 
are equivalent to certain special cases of 
the \IC for the polynomial algebra $\cxz$ and 
the commuting differential operators $\Theta$ 
defined above. The main results of this subsection 
are Theorems \ref{VC2IC}, \ref{JC2IC} and \ref{AJC2IC}.

Finally, in Subsection \ref{S3.3} we give a different 
description for the image $\im\Theta$ of 
the differential operators $\Theta$ above. 
The main result of this subsection is Theorem 
\ref{2nd-description} for which we give 
two proofs. From the second proof, we will  
see that the \IC is actually also related 
with the multidimensional Laplace transformations 
of polynomials (See Corollary \ref{L2nd-description} and 
Conjecture \ref{RestateIC}). \\

{\bf Acknowledgment} The author would like to thank Professor Arno van den Essen and the anonymous referee for valuable suggestions to improve the earlier version of this paper.

\renewcommand{\theequation}{\thesection.\arabic{equation}}
\renewcommand{\therema}{\thesection.\arabic{rema}}
\setcounter{equation}{0}
\setcounter{rema}{0}

\section{\bf The Image Conjecture of Commuting Differential Operators of Order One with Constant Leading Coefficients}
\label{S2}

The most interesting case of the {\it image} conjecture ({\bf IC}), Conjecture \ref{IC}, probably is the case when the commutative algebra $\cA$ is a field $K$ of characteristics zero. For example, as we will show later in Subsection \ref{S3.2}, the Jacobian conjecture is actually equivalent to a very special case of the \IC with $\cA=\bC$. But, in order to get a good reduction for the {\bf IC}, we need to consider a slightly more general case. Namely, in this section we consider the case of the \IC under the following two assumptions.

\begin{enumerate} 
\item[$({\bf C_1})$] \label{2Cs} the commutative algebra $\cA$ is a commutative algebra over a field $K$ of 
characteristics zero;

\item[$({\bf C_2})$] the commuting differential operators in 
 $\cC$ have the form 
$h(z)+\sum_{i=1}^n c_i\p_i$ with 
$h(z)\in \cAz$ and the leading coefficients 
$c_i\in K$ (not just in $\cA$) for any $1\le i\le n$.
\end{enumerate}

Throughout this section, $K$ always denotes a field of characteristic zero and $\cA$ a commutative algebra over $K$. 
We denote by $\bD_K[z]$ the set of differential operators 
of $\cAz$ of the form 
$h(z)+\sum_{i=1}^n c_i\p_i$, where  
$h(z)\in \cAz$ and the leading coefficients 
$c_i\in K$ $(1\le i\le n)$. Furthermore, 
all the terminologies and notations introduced 
in the Introduction will also be in 
force throughout this paper. 

We start with the following simple lemma.

\begin{lemma}\label{L2.1.1}
Let $1\le k\le n$ and $\Phi_i=\p_i-h_i(z)$ $(1\leq i \leq k)$ 
with $h_i(z)\in \cAz$. Then  
the differential operators $\Phi_i$ $(1\le i\le k)$ commute with one another iff there exists a $q(z)\in \cAz$ such that 
$h_i=\p_i (q)$ for all $1\le i\le k$.
\end{lemma}

\pf First, for any $1\le i, j\le k$, we have the following 
identity for the commutator of $\Phi_i$ and $\Phi_j$:  
\begin{align*}
[\p_i-h_i, \p_j-h_j]=\p_j(h_i)-\p_i(h_j).
\end{align*}
Therefore, the differential operators $\Phi_i$'s 
commute with one another iff, for any $1\le i, j\le k$, 
\begin{align}\label{L2.1.1-pe1}
\p_j(h_i)=\p_i(h_j). 
\end{align}

Since $\cA$ is an $K$-algebra, hence also an $\bQ$-algebra, by Poincar\'e's Lemma, we know that Eq.\,(\ref{L2.1.1-pe1}) holds for all $1\le i, j\le k$ iff there exists a $q(z)\in K[z]$ such that $h_i(z)=\p_i(q)$ for all $1\le i\le k$.  \epfv

Note that $\bD_K[z]$ is a $K$-vector space and,  
for any subset $\cC \subset \bD_K[z]$, 
it is easy to check 
that $\im\cC=\im V$, where $V$ is the 
$K$-subspace of $\bD_K[z]$ spanned by elements of $\cC$ 
over $K$. Therefore, without losing any generality, we may freely replace the commuting subsets in the \IC by 
commuting $K$-subspaces of $\bD_K[z]$.

\begin{lemma}\label{L2.1.2}
Let $V$ be a commuting $K$-subspace of $\bD_K[z]$ that contains at least one $($nonzero$)$ differential operator of order one. Then up to an automorphism of $\cAz$, 
there exist $1\le k\le n$,  
$q(z)\in \cAz$, and finitely many polynomials 
$g_i(z)\in \cAz$ $(i\in I)$ such that  
\begin{enumerate}
\item $g_i(z)$ $(i\in I)$ does not involve  
$z_j$ for any $1\le j\le k$, i.e.\@ 
$\p_j g_i(z)=0$ for all $i\in I$ and $1\le j\le k$;
\item $\im V$ is the same as the image of the collections of 
the differential operators $\p_j-\p_j(q(z))$ $(1\le j\le k)$ and the multiplication operators by $g_i(z)$ $(i\in I)$.
\end{enumerate}
\end{lemma}

\pf First, for any $u=(a_1, a_2, ... , a_n)\in K^n$, set 
$u\p\!:=\sum_{i=1}^n a_i\p_i$. 
Let $U=\{ u\in K^n \,|\, u\p+h(z)\in V 
\mbox{ for some $h(z)\in \cAz$} \}$ and 
$V_0$ the set of all differential operators in $V$ of order zero, i.e.\@  the multiplication operators 
by elements of $\cAz$. 
Then by the fact that $V$  
is a vector space over $K$ 
(actually an $K$-subspace of $\bD_K[z]$), 
it is easy to check that both $U$ and $V_0$ 
are also vector spaces over $K$. 

Furthermore, by the assumption on $V$ in the lemma, 
we have $U\ne 0$. Let $u_1, u_2, ..., u_k$ 
be a basis of $U$. Then  
up to a changing of coordinates, we may assume that 
$u_i=e_i$ $(1\le i\le k)$, where the $e_i$'s are the standard 
basis vectors of $K^n$. Under this assumption, 
we have $\p_i-h_i(z) \in V$ $(1\le i\le k)$ 
for some $h_i(z)\in \cAz$. 
Since $V$ is a commuting subset of $\bD_K[z]$,  
by Lemma \ref{L2.1.1} there exists a $q(z)\in \cAz$ 
such that $h_i=\p_i (q)$ for all $1\le i\le k$. 
We denote by $V_1$ the $K$-subspace of 
$V$ spanned by $\p_i-\p_i(q)$ 
$(1\le i\le k)$ over $K$. 

Next we consider $V_0$, which can be viewed as 
a $K$-subspace of $\cAz$. Let $g_i(z)\in \cAz$ 
$(i\in J)$ be a $K$-basis of $V_0$. 
Then the polynomials $g_i(z)\in \cAz$ 
$(i\in J)$ also generate the ideal 
$V_0K[z]$ of $K[z]$. On the other hand, 
by Hilbert's theorem, we know that 
$K[z]$ is Noetherian. Therefore, there exists 
a finite subset $I\subset J$ such that 
$g_i(z)\in \cAz$ $(i\in I)$ also 
generate the ideal $V_0 K[z]$.

Now, by the fact that $\im V=\im V_1+\im V_0=\im V_1+V_0 K[z]$, 
it is easy to check that $\im V$ is the same as the image of the collections of the differential operators $\p_j-\p_j(q(z))$ $(1\le j\le k)$ and the multiplication operators by $g_i(z)$ $(i\in I)$.

Finally, for any $i\in I$ and $1\le j\le k$, 
since $V$ is commuting, we have 
$[\p_j-h_j, g_i(z)]=\p_j(g_i)=0$. Therefore, 
$g_i(z)$ $(i\in I)$ does not involve  
$z_j$ for any $1\le j\le k$. 
\epfv   

Note that, if a commuting $K$-subspace 
$V\subset \bD_K[z]$ does not contain 
any nonzero differential operators 
of order one, i.e.\@  $V$ contains only multiplication operators by elements of $K[z]$, we can identify $V$ with 
a $K$-subspace of $K[z]$. 
Then it is easy to see that $\im V$ is the same as 
the ideal $VK[z]$ of $K[z]$ generated by elements of $V$. 
In this case, the \IC holds trivially.
 
Therefore, we may assume that $V$ contains 
at least one nonzero differential operators 
of order one. Furthermore, by Lemma \ref{L2.1.2}, 
without changing the image, we may also assume that 
$V$ is linearly spanned over $K$ by the differential 
operators $\p_j-\p_j(q(z))$ $(1\le j\le k)$ 
for a $q(z)\in K[z]$ and 
finitely many multiplication operators by 
$g_i(z)\in \cAz$ $(i\in I)$ such that  
the polynomials $g_i(z)$ $(i\in I)$ do  
not involve the variables  
$z_j$ $(1\le j\le k)$.

Throughout the rest of this section, 
we will fix a commuting subspace $V\subset \bD_K[z]$ 
as above. Set $z'\!:=(z_1, ..., z_k)$ and 
$z''\!:=(z_{k+1}, z_{k+2}, ..., z_n)$. 
Then we have $g_i(z)\in \cA[z'']$ $(i\in I)$. 
Let  $\mathcal I_1$ (resp.\,$\mathcal I_2$) 
be the ideal of $\cA[z'']$ (resp.\,$\cAz$) generated by 
$g_i(z)$ $(i\in I)$. Set $\cB\!:=\cA[z'']/\mathcal I_1$ 
and denote by $\pi_1: A[z'']\to \cB$ the quotient map. 

Since $g_i(z)\in \cA[z'']$ $(i\in I)$, the quotient space $\cAz/\mathcal I_2$ can (and will) be naturally identified 
with $\cB[z']$, and the quotient map  
\begin{align*}
\pi_2: \cAz=\cA[z''][z'] \to \cAz/\mathcal I_2 \simeq 
\cB[z']
\end{align*} 
can be viewed as the linear extension of the quotient map 
$\pi_1: \cA[z''] \to \cB$ to the polynomial algebras in $z'$ over 
$\cA[z']$ and $\cB$. 

For convenience, we introduce the following shorter notations. 
For any $u(z)\in \cAz$, we 
denote by $\bar u(z)$ the image of $u(z)$ under the quotient map 
$\pi_2$, i.e.\@  $\bar u(z)\!:=\pi_2(u(z))$. 
For any $1\le i\le k$, we set $h_i(z)\!:=\p_i(q(z))$; 
$\Phi_i\!:=\p_i-h_i(z)$ and $\Psi_i\!:=\p_i-\bar h_i(z)$.
Note that $\Psi_i$ $(1\le i\le k)$ are commuting differential operators of $\cB[z']$ of order one with leading coefficients lying in the base field $K$. Finally, set 
$\overline{V}\!:= \mbox{Span}_K \{ \Psi_j \,|\, 1\le j\le k\}$. 

\begin{lemma}\label{L2.3a}
With the notations fixed as above, we have

$(a)$ for any $1\le j\le k$ and $f(z)\in \cAz$, 
\begin{align}
\pi_2 (\Phi_j f)&=\Psi_j \bar f. \label{L2.3a-e1}
\end{align}

$(b)$ 
\begin{align}\label{E4-pi}
\im V=\pi_2^{-1}(\im \overline{V}).
\end{align}
\end{lemma}
\pf $(a)$ Note first that, since $g_i(z)$ $(i\in I)$ are independent on $z'$, the ideal $\mathcal I_2\subset \cA[z]$ is preserved by the derivations $\p_j$ $(1\le j\le k)$. Consequently, for any $1\le j\le k$, we have 
\begin{align}\label{L2.3a-pe1}
\pi_2\circ\p_j=\p_j\circ \pi_2.
\end{align}

For any $1\le j\le k$, by Eq.\,(\ref{L2.3a-pe1}) and the fact that $\pi_2$ is a homomorphism of algebras, we have 
\begin{align*}
\pi_2 (\Phi_j f)&=\pi_2 (\p_j f-h_j f)=\pi_2 (\p_j f)-\pi_2(h_j f)\\
&=\p_j\pi_2(f)-\pi_2(h_j)\pi_2(f)= \p_j \bar f-\bar h_j \bar f\\
&=(\p_j-\bar h_j)\bar f=\Psi_j \bar f. 
\end{align*}
Therefore, we have $(a)$.
 
To show $(b)$, let $f\in \im V$ and write it as 
$f=\sum_{j=1}^n\Phi_j a_j + \sum_{i\in I} g_i b_i $ 
for some $a_j , b_i \in \cAz$ $(i\in I)$. 
Then by Eq.\,(\ref{L2.3a-e1}) and the fact that $\pi_2(g_i)=0$ for any $i\in I$, we have 
$\pi_2(f)=\sum_{j=1}^n \Psi_j \bar a_j \in \im \overline{V}$ and
$f\in  \pi_2^{-1}(\im\overline{V})$.

Conversely, let $f\in \pi_2^{-1}(\im\overline{V})$, i.e.\@  
$\bar f\in \im\overline{V}$. Since $\pi_2$ is surjective, we may write $\bar f=\sum_{j=1}^n \Psi_j \bar a_j$ for some 
$a_j\in \cAz$ $(1\le j\le k)$. Then  
by Eq.\,(\ref{L2.3a-e1}), we have 
$\bar f=\sum_{j=1}^n \pi_2( \Phi_j  a_j)$. 
Hence  
$f=\sum_{j=1}^n \Phi_j a_j(z)+u$ for some 
$u\in {\mathcal I}_2=\Ker \pi_2$. 
Since ${\mathcal I}_2$ is the ideal of $\cAz$ generated 
by $g_i(z)$ $(i\in I)$, we have 
${\mathcal I}_2\subset \im V$. Therefore, 
$u\in \im V$. Hence we also have $f\in  \im V$. 
\epfv

\begin{lemma}\label{L2.1.3}
With $V$ and related notations fixed as above, we have that  
the \IC holds for $\cAz$ and the commuting differential operators $V$, iff the \IC holds for $\cB[z']$ and the commuting differential operators 
$\{\Psi_j\,|\, 1\le j\le k\}$. 
\end{lemma}

\pf $(\Rightarrow)$ Assume that the 
\IC holds for $\cAz$ and the commuting differential 
operators $V$. Take any $\bar f, \bar g \in \cB[z']$ 
with $\bar f^m \in \im \{\Psi_j \, |\, 1\le j\le k\}
=\im \overline{V}$ for all $m\ge 1$. 
By Eq.\,(\ref{E4-pi}), we have $f^m \in \im V$ for all 
$m\ge 1$. Therefore, we have $f^m g \in \im V$ 
when $m\gg 0$. By Eq.\,(\ref{E4-pi}) again, $\bar f^m \bar g=\pi_2(f^mg)\in \im \overline{V}=\im \{\Psi_j \, |\, 1\le j\le k\}$ 
when $m\gg 0$.

$(\Leftarrow)$ Assume that the 
\IC holds for $\cB[z']$ and the commuting differential 
operators $\{\Psi_j \, |\, 1\le j\le k\}$, 
hence also for $\overline{V}$. For any $f, g\in \cAz$ with 
$f^m\in \im V$ for all $m\ge 1$, by Eq.\,(\ref{E4-pi}) 
we have $\bar f^m \in \im \overline{V}$ for any $m\ge 1$. 
Therefore, $\pi_2(f^mg)=\bar f^m \bar g  \in \im \overline{V}$ 
when $m\gg 0$. By Eq.\,(\ref{E4-pi}) again, 
we have $f^mg\in \im V$ when $m\gg 0$. 
\epfv

By Lemmas \ref{L2.1.1}--\ref{L2.1.3}, it is easy to see that 
we have the following reduction on the \IC under the conditions 
$({\bf C_1})$ and $({\bf C_2})$ on page $\pageref{2Cs}$. 

\begin{theo}\label{T2.1.4}
To prove or disprove the \IC under the conditions 
$({\bf C_1})$ and $({\bf C_2})$, it is enough to consider 
the \IC for polynomial algebras $\cAz$ $($in $n$ variables$)$
for the commuting differential operators 
$\cC=\{\p_i-\p_i(q(z))\,|\, 1\le i\le n\}$ 
with $q(z)\in \cAz$.
\end{theo}

Two remarks on the \IC are as follows. First, it does not hold in general for commuting differential operators of order one with non-constant leading coefficients. 

\begin{exam}\label{E2.1.5}
Consider the polynomial algebra $\bC[t]$ in one variable $t$. Let $\Phi\!:=td/ dt -1$.  Then it is easy to check that
for any $f(t) \in \bC[t]$ $f(t)\in \im \Phi$ iff 
$f'(0)=0$, i.e.\@  $f(t)$ has no degree-one term.

Now let $f(t)=1+t^2$. Then $f^m(t)\in \im \Phi$ 
for each $m\ge 1$. But $t f^m(t) \not \in \im \Phi$ for each
$m\ge 1$ since 
the degree-one terms of $t f^m(t)$ are always $t$. 
Therefore, $\im \Phi$ is not a Mathieu subspace of $\bC[t]$ 
and the \IC fails in this case.
\end{exam}

Second, the \IC is also false in general when the base algebra $\cA$ is an algebra over a field of characteristic $p>0$.

\begin{exam}\label{E2.1.6}
Let $K$ be any field of characteristic $p>0$ and $x$ a free variable. Consider the differential operator 
$\Lambda\!:=d/ dx$ of the polynomial algebra $K[x]$.  
Then it is easy to check that $\im \Lambda$ is linearly spanned over $K$ by the monomials $x^m$ $($$m\in \bN$ and 
$m\not \equiv -1$ $($\mbox{\rm mod} $p))$. In particular, 
the constant polynomial $1\in \im \Lambda$. Thus by taking 
$f=1$ and $g=x^{p-1}$, we see that 
the \IC fails for the positive 
characteristic case.
\end{exam}

Next we further assume $\cA=K$ and consider 
some $\mcD$-module structures on $\im\cC'$ and 
$K[z]/\im\cC'$ for any subset $\cC'$ of the commuting 
differential operators $\cC$ in Theorem \ref{T2.1.4}.

We believe that most of the results to be discussed below,  
if not all, are known in the theory of $\mcD$-modules. But, for the sake of completeness and also due to lack of more specific references in the literature, we give a more detailed account here. It is a little surprising to see that the problem raised by the \IC has not gotten any attention 
from the point of view of $\mcD$-modules. 

Note first that up to a permutation of the variables $z_i$  
we may assume $\cC'=\{\Phi_i\,|\, 1\le i\le k\}$ 
for some $1\le k\le n$. For the simplicity of notation, 
throughout the rest of this subsection, 
we fix an $1\le k\le n$ and still 
denote $\cC'$ by $\cC$.

Let $\mcD[z]$ be the Weyl algebra of the polynomial 
algebra $\Kz$, i.e.\@ the subalgebra 
of $\mbox{\rm End}_K(\Kz)$ (the algebra of all  
$K$-linear maps from $\Kz$ to $\Kz$)   
generated by the derivations 
$\p_i$ $(1\le i\le n)$ and 
the multiplication operators 
by $z_i$ $(1\le i\le n)$.  
We first define an action of $\mcD[z]$ on 
$\Kz$ by setting 
\begin{align}
z_i\ast f&\!:=z_i f, \label{s-act1}\\
\p_i\ast f&\!:=\Phi_i f=\p_i f-f\p_i q 
\label{s-act2}
\end{align}
for all $1\le i\le n$ and $f\in \Kz$.

Since $\Phi_i$ $(1\le i\le n)$ commute with one another, 
it is easy to see that $K[z]$ with the actions defined above 
forms a module of the Weyl algebra $\mcD[z]$. 
We denote this $\mcD[z]$-module by $\cM$.  

The $\mcD[z]$-module $\cM$ can also be constructed as follows. 
First let $\cM_q\!:=K[z]e^{-q(z)}$. We may formally view 
elements of $K[z]e^{-q(z)}$ as formal power series in $z$ 
over $K$.  Then the standard action of $\mcD[z]$ 
on the formal power series algebra $K[[z]]$ induces 
an action of $\mcD[z]$ on $\cM_q$. More precisely, 
for any $1\le i\le n$ and $f(z)\in K[z]$, we have
\begin{align}
z_i\cdot (f(z)e^{-q(z)})&= (z_i f(z))e^{-q(z)}, \label{q-act1}\\
\p_i\cdot (f(z)e^{-q(z)})&=e^{-q(z)}(\p_i f-f\p_i q)
=e^{-q(z)}(\Phi_i\, f).\label{q-act2}
\end{align}

Since $\cM_q$ is obviously closed under the actions above, 
it is also a $\mcD[z]$-module. Furthermore, 
let $E: \cM \to \cM_q$ be the $K$-linear map that maps 
any $f(z)\in K[z]$ to $f(z)e^{-q(z)} \in \cM_q$. Then   
by Eqs.\,(\ref{s-act1})--(\ref{q-act2}),
it is straightforward to check that the following lemma holds.

\begin{lemma}\label{L2.1.6}
\mbox{} 
$(a)$ The map $E$ is an 
isomorphism of $\mcD[z]$-modules.

$(b)$ For any $1\le i\le n$, we have
\begin{align}\label{L2.1.6-e1}
E(\p_i\ast \cM)= \p_i\cdot \cM_q,
\end{align}
where $\p_i\ast \cM$ and $\p_i\cdot \cM_q$ are the images of 
$\cM$ and $\cM_q$ under the actions of 
$\p_i\ast$ and $\p_i\cdot$, respectively.

$(c)$ Let $\cC=\{\Phi_i\,|\, 1\le i\le k\}$ as fixed 
before. Then viewing $\im\cC$ as a subspace of $\cM$, 
we have $\im\cC=\sum_{i=1}^k \p_i\ast \cM$.
\end{lemma}
   
\begin{lemma}\label{L2.1.7}
$\cM$ is a holonomic module of the Weyl algebra 
$\mcD[z]$ with multiplicity 
$e(\cM) \le \max\{1, d^n\}$, 
where $d\!:=\deg q(z)$. 
\end{lemma}

Since $\cM\simeq \cM_q$ as $\mcD[z]$-modules, it is enough to show the lemma for $\cM_q$. First, it is easy to just show that $\cM_q$ is holonomic as follows. \\
 
\underline{\it First proof:} Let $I_q$ be the left ideal of 
$\mcD[z]$ consisting of differential operators that annihilate 
$e^{-q(z)}$. Then it is easy to see that 
$\cM_q\simeq \mcD[z]/I_q$ as $\mcD[z]$-modules. 
Since $e^{-q(z)}$ is an {\it $($hyper$)$exponential} 
function, by Theorem $2.6$ in \cite{C} (p.\,183) 
we obtain that $e^{-q(z)}$ is a holonomic function, 
i.e.\@ $\mcD[z]/I_q$ is a holonomic $\mcD[z]$-module. 
Hence so is $\cM_q$.
\epfv

A more straightforward proof of Lemma \ref{L2.1.6} which also gives the information on the multiplicity $e(\cM)$ can be given as follows. \\

\underline{\it Second Proof:} First, if $d=0$, then 
$q(z)$ is a constant and  
$\cM_q$ is isomorphic to the standard 
$\mcD[z]$-module $K[z]$. It is well-known that $K[z]$ 
is holonomic with multiplicity $e(K[z])=1$. Hence
the lemma holds in this case. 

Assume $d\ge 1$ and introduce a filtration on $\cM_q$ by setting,  
for any $m\ge 0$,
\begin{align}
\Gamma_m\!:= \{ u(z)e^{-q(z)} \, | \, \deg u(z) \le md \} 
\end{align}  

First, it is easy to see that 
$\cup_{m\ge 0} \Gamma_m= K[z]e^{-q(z)}=\cM_q$ since 
$md\to \infty$ as $m\to \infty$.

Second, for any $m\ge 0$, 
$f(z)e^{-q(z)}\in \Gamma_m$ and $1 \le i \le n$, by 
Eqs.\,(\ref{q-act1})-(\ref{q-act2}),  
we have $z_i\cdot(f(z)e^{-q(z)})=(z_i f(z))e^{-q(z)}
\in \Gamma_{m+1}$ since
\begin{align*}
\deg(z_i f(z))
=\deg f+1\le md + 1 \le (m+1)d,
\end{align*} 
and
\begin{align*}
\p_{i}\cdot (f(z)e^{-q(z)})=
e^{-q(z)}(\p_i f-f\p_i q) \in \Gamma_{m+1} 
\end{align*}
since $\deg\, (\p_i f-f\p_i q) \le \deg f+d -1\le md+d-1 
<d(m+1)$.

Third, for any fixed $m\ge 0$, 
$\dim \Gamma_m$ is the same as the dimension of 
the subspace of $\cz$ consisting of the polynomials of 
degree less or equal to $md$. Therefore, for any 
$r\in \bR$ such that 
$r\ge \frac 1{n!}\binom nk n^{n-k}$ 
for all $0\le k\le n-1$, we have 
\allowdisplaybreaks{
\begin{align}
\dim \Gamma_m&=\binom{n+md}{n}\le \frac{(n+md)^n}{n!} \\
&=
\frac{(md)^n}{n!}+\frac 1{n!}
\sum_{k=0}^{n-1} \binom nk n^{n-k}(md)^k \nno \\
&\le \frac{(md)^n}{n!}+ rn (md)^{n-1}
\le \frac{d^n m^n}{n!}+c(m+1)^{n-1}, \nno 
\end{align} }
where $c =rnd^{n-1}$. 

Then by Theorem $5.4$ (p.\,12) in \cite{B} or 
Lemma $3.1$ (p.\,91) in \cite{C}, 
$\cM_q$ is a holonomic $\mcD_K$-module with 
multiplicity $e(\cM_q)\le d^n$.    
\epfv

Now let $1\le k\le n$ and $\cC$ as fixed before. 
Set $z'\!:=(z_1, ..., z_k)$ and $z''\!:=(z_{k+1}, ..., z_n)$. 
Denote by $\mcD[z'']$ the Weyl algebra of the polynomial algebra
$K[z'']$. Note that the $\mcD[z]$-module $\cM$ is also a 
$\mcD[z'']$-module since $\mcD[z'']$ is a subalgebra of the Weyl algebra $\mcD[z]$. 

\begin{theo}\label{T2.1.8}
$(a)$ $\im\cC$ as a subspace of $\cM$ is a  
$\mcD[z'']$-submodule of the 
$\mcD[z'']$-module $\cM$.

$(b)$ The quotient 
$\cM/\im\cC$ is a holonomic $\mcD[z'']$-module. 
\end{theo}

\pf $(a)$ Note that for any $k+1\le j\le n$, 
the differential operator $\Phi_j$ and 
the multiplication operator by $z_j$ 
commute with $\Phi_i$ $(1\le i\le k)$. Therefore, 
$\im\cC$ is closed under the actions of 
$\Phi_j$ and $z_j$ $(k+1\le j\le n)$ and 
forms a $\mcD[z'']$-submodule of $\cM$.

$(b)$ By Lemma \ref{L2.1.7}, we know that  
$\cM$ is a holonomic $\mcD[z]$-module. Then it follows 
from Theorem $6.2$ (p.\,16) in \cite{B} that 
the quotient $\cM/(\sum_{i=1}^k \p_i\ast\cM)$ 
is a holonomic $\mcD[z'']$-module. 
Combining this fact with Lemma \ref{L2.1.6}, 
$(c)$, we see that the statement $(b)$ 
in the theorem holds. \epfv

\begin{corol}\label{C2.1.9} 
When $k=n$, $\im\cC$ is a finite co-dimensional 
$K$-subspace of $K[z]$, i.e.\@  the $K$-vector space 
$K[z]/\im\cC$ is of finite dimension.   
\end{corol}

\pf Note that when $k=n$, we have $\mcD[z'']=K$. 
Then by Theorem \ref{T2.1.8}, $(b)$, we know that  
$K[z]/\im\cC$ is a holonomic $K$-module. Hence,  
$K[z]/\im\cC$ must be finite dimensional 
over $K$. 
\epfv

\renewcommand{\theequation}{\thesection.\arabic{equation}}
\renewcommand{\therema}{\thesection.\arabic{rema}}
\setcounter{equation}{0}
\setcounter{rema}{0}

\section{\bf The Vanishing Conjecture and the Jacobian Conjecture in Terms of the Image Conjecture} \label{S3}

In this section, we study the relations of the 
{\it vanishing} conjecture ({\bf VC}), 
Conjecture \ref{VC}, and  the {\it Jacobian} 
conjecture ({\bf JC}) with the {\it image} 
conjecture ({\bf IC}), Conjecture \ref{IC}. 
But first, we need to fix more notations. 

Let $\xi=(\xi_1, \xi_2, ..., \xi_n)$ be another 
$n$ commutative free variables which also commute 
with the free variables $z=(z_1, z_2, ..., z_n)$. 
For any $1\le i\le n$, set $\Theta_i\!:=\xi_i-\p_i$ and 
$\Theta\!:=\{\Theta_i\,|\, 1\le i\le n\}$.
Note that $\Theta$ is a commuting set of 
differential operators of order one  
with constant leading coefficients of 
the polynomial algebra $\cxz$. 

In Subsection \ref{S3.2}, we show  
that the \VC and the \JC are equivalent 
to some special cases of the \IC for the polynomial algebra $\cxz$ and the commuting differential operators $\Theta$. 
The main results are given in 
Theorems \ref{VC2IC}, \ref{JC2IC} and \ref{AJC2IC}.   
 
One crucial result needed for the proofs of the theorems above 
is Theorem \ref{Main-Thm-1}. So we devote Subsection \ref{S3.1} to prove a stronger version, Theorem \ref{D-Taylor}, of this result. We believe that these two results are interesting and important on their own rights, so we have formulated them as theorems. 

In Subsection \ref{S3.3} we give another description 
for $\im \Theta$ which does not involve differential operators 
(See Theorem \ref{2nd-description}). From the second proof of the theorem, we will see that the $\IC$ is also connected with integrals of polynomials such as the multidimensional Laplace transformations.

\subsection{Images of Commuting Differential Operators $\Theta$}
\label{S3.1}

First, let us fix the following conventions and notations that will be used throughout the rest of this paper.
 
\begin{enumerate}
  \item For any $n\ge 1$, we define a partial order on $\bN^n$ by setting, for any $\alpha, \beta\in \bN^n$, 
   $\alpha \ge \beta$ if, for each $1\le i\le n$, the $i^{\rm th}$ component of $\alpha$ is greater or equal to the $i^{\rm th}$ component of $\beta$.
  \item We will freely use some multi-index notations. For instance, for any $\alpha=(k_1, k_2, ..., k_n)\in \bN^n$, we set
 \begin{align*}
\alpha!&=k_1!k_2!\cdots k_n!.\\
z^\alpha&=z_1^{k_1}z_2^{k_2}\cdots z_n^{k_n}.\\
\p^\alpha&=\p_1^{k_1}\p_2^{k_2}\cdots \p_n^{k_n}.
\end{align*} 
\end{enumerate}

Next, we define a linear map $\mcE: \cxz\to \cz$ by setting,   
for any $g(\xi) \in \bC[\xi]$ and $h(z) \in \bC[z]$,  
\begin{align}\label{NewDef-E}
\mcE\big(g(\xi)h(z)\big)\!:= g(\p)h(z).
\end{align}

In other words, for any $f(z, \xi)\in \bC[z, \xi]$, 
$\mcE (f(z, \xi))\in \bC[z]$ is obtained 
by putting the free variables $\xi_i$'s on the most left 
in each monomial of $f(z, \xi)$ and then replacing 
$\xi_i$ by $\p_i$ and applying them to $z_i$'s on the right. For example, we have
{\it 
\begin{enumerate}
  \item $\cE(1)=1;$
\item $\mcE(z_1^n \xi_1^2)=\p_1^2(z_1^n)=n(n-1)z_1^{n-2}$ 
for each $n\ge 2$; 
\item $\mcE(a(z))=a(z)$ for each  $a(z)\in \cz$; 
\item $\mcE(\xi^\alpha)= \p^\alpha(1)=0$  for each  
$0\neq\alpha\in \bN^n$. 
\end{enumerate}  
}

Next, we are going to derive another description of 
the subspace $\Ker \mcE\subset$$\cxz$ by using 
the differential operators 
$\Theta=\{\Theta_i=\xi_i-\p_i\,|\, 1\le i\le n\}$ 
(as fixed in the first paragraph of this section) 
of the polynomial algebra $\cxz$. But, for convenience, 
in the rest of this paper, we also use $\Theta$ to denote 
the ordered $n$-tuple $(\Theta_1, \Theta_2, ..., \Theta_n)$, 
so the notations such as $\Theta^\alpha$ $(\alpha\in \bN^n)$ 
will make perfect sense.

\begin{theo}\label{Main-Thm-1}
$\Ker \mcE=\im \Theta$.
\end{theo}

Actually, a much more explicit theorem, 
Theorem \ref{D-Taylor}, can be formulated and proved. 
Theorem \ref{Main-Thm-1} follows immediately from 
Theorem \ref{D-Taylor} and lemma \ref{Main-Thm-1-L} 
below. Indeed, assume $f \in \Ker \mcE$. Then by Theorem
\ref{D-Taylor}, $a_0(z) = 0$ and thus $f \in \im \Theta$.

\begin{theo}\label{D-Taylor}
$(1)$ Any $f(\xi, z)\in \cxz$ can be written uniquely as a sum of the form 
\begin{align}\label{D-Taylor-e1}
f(\xi, z)=\sum_{\alpha\in \bN^n} \frac 1{\alpha!} \Theta^\alpha a_\alpha(z)
=\sum_{\alpha\in \bN^n} \frac 1{\alpha!} (\xi-\p )^\alpha a_\alpha(z)
\end{align}
for some $a_\alpha(z)\in \bC[z]$ $(\alpha\in \bN^n)$.

$(2)$ With the notation above, for any 
$\alpha=(\alpha_1, \alpha_2, ..., \alpha_n)\in \bN^n$,
\begin{align}\label{D-Taylor-e2}
a_\alpha(z)=\mcE(\p_\xi^\alpha f),
\end{align}
where $\p_\xi^\alpha=\prod_{i=1}^n \p_{\xi_i}^{\alpha_i}$.

In particular, we have $a_0(z)=\mcE(f)$.
\end{theo}

To prove Theorem \ref{D-Taylor}, we first need to show the following lemma which is actually the easy 
part of Theorem \ref{Main-Thm-1}.

\begin{lemma}\label{Main-Thm-1-L}
$ \im\Theta\subseteq \Ker\mcE$.
\end{lemma}

\pf 
For any $f(\xi, z)\in \im\Theta$, we need to show $\mcE(f)=0$.
By the linearity, we may assume 
$f(\xi, z)=\Theta_i (a(\xi)b(z))$ for some $1\le i\le n$, 
$a(\xi)\in \bC[\xi]$ and $b(z)\in \cz$. Then 
by Eq.\,(\ref{NewDef-E}), we have
\allowdisplaybreaks{
\begin{align*}
\mcE(\Theta_i(a(\xi)b(z)) )&=\mcE(\, \xi_i a(\xi)b(z)-\p_i(a(\xi)b(z))\, )\\
&=\mcE( a(\xi)\xi_i b(z)- a(\xi)\p_i b(z) ) \\
&=a(\p)\p_i b(z)- a(\p)\p_i b(z)=0.
\end{align*} }
Therefore, we have $f(\xi, z)\in \Ker\mcE$ and 
$\im\Theta\subseteq \Ker\mcE$.
\epfv

\underline{\it Proof of Theorem \ref{D-Taylor}}: 
First, we show that Eq.\,(\ref{D-Taylor-e1}) does hold 
for some $a_\alpha(z)\in \bC[z]$ $(\alpha\in \bN^n)$. 
By the linearity, we may assume that 
$f(\xi, z)=\xi^\beta h(z)$ for 
some $\beta\in \bN^n$ and 
$h(z)\in \bC[z]$.

Consider
\begin{align}
\xi^\beta h(z)&=(\xi-\p +\p )^\beta h(z)
=(\Theta+\p )^\beta h(z)  \\
&=\sum_{\substack{\alpha \in \bN^n \\ \alpha \le \beta}} 
\frac {\beta!}{\alpha!(\beta-\alpha)!} \Theta^\alpha 
\p ^{\beta-\alpha} h(z), \nno
\end{align}
which is of the desired form of 
Eq.\,(\ref{D-Taylor-e1}) for $f(\xi, z)=\xi^\beta h(z)$  
with 
\begin{align*}
a_\alpha(z) =
\begin{cases}
\frac {\beta!}{(\beta-\alpha)!}  \, 
\p ^{\beta-\alpha} h(z) &\mbox{ if } \alpha\le \beta; \\
0 &\mbox{ otherwise.}
\end{cases}
\end{align*}

Next, we show the uniqueness of 
Eq.\,(\ref{D-Taylor-e1}) as follows.

Note first that for any $\alpha,\beta \in \bN^n$, we have 
\begin{align*}
\p_\xi^\beta \left(\frac 1{\alpha!} \Theta^\alpha a_\alpha(z)\right)
=
\begin{cases}
\frac 1{(\alpha-\beta)!}\, \Theta^{\alpha-\beta} a_\alpha(z)
& \mbox{ if } \beta\le \alpha; \\
0 &\mbox{ otherwise}.
\end{cases}
\end{align*}

Then by applying $\p_\xi^\beta$ 
to Eq.\,(\ref{D-Taylor-e1}), we get 
\begin{align}
\p_\xi^\beta f(\xi, z)=\sum_{\gamma \in \bN^n} 
\frac 1{\gamma!} \Theta^\gamma(a_{\gamma+\beta}(z)).
\end{align}

Applying $\mcE$ to the equation above, by Lemma 
\ref{Main-Thm-1-L} we see that 
$a_\beta(z)=\mcE(\p_\xi^\beta f(\xi, z))$ 
for each $\beta\in \bN^n$, which is exactly 
Eq.\,(\ref{D-Taylor-e2}) with the index $\alpha$ 
replaced by $\beta$. Hence, the theorem follows. 
\epfv

\begin{rmk}
A more conceptual proof of Theorem \ref{D-Taylor} will be given in \cite{T-Deform}. Eq.\,$(\ref{D-Taylor-e1})$ actually corresponds to the Taylor series in a deformation of the polynomial algebra $\cxz$. It can also be viewed as a twisted version of the usual Taylor series of $f(\xi, z)\in \cxz$ by the commuting differential operators $\Theta$.
\end{rmk}

\subsection{The Vanishing Conjecture and the Jacobian Conjecture in Terms of the Image Conjecture}\label{S3.2}

Now we are ready to show that some relations of the \VC and 
the \JC with the {\bf IC}. First, the relation between the \VC and 
the \IC is given as follows.

\begin{theo}\label{VC2IC}
For any $\Lambda(\xi)\in \bC[\xi]$ and $P(z)\in \bC[z]$,
the following two statements are equivalent:

$(a)$ the \VC holds for $\Lambda=\Lambda(\p)$ 
and $P(z), g(z)\in \cz$;

$(b)$  the \IC holds for $f(\xi, z)\!:=\Lambda(\xi)P(z)$ and 
$g(z)$ in $\cxz$.
\end{theo} 

\pf For any $m\ge 1$, by Eq.\,(\ref{NewDef-E}) and 
Theorem \ref{Main-Thm-1}, 
we have that  
$\Lambda(\p)^m(P^m(z))=0$ iff 
$\mcE(\Lambda^m(\xi)P^m(z))=\mcE(f^m(\xi, z))=0$ 
iff $f^m(\xi, z)\in \Ker (\mcE)=\im\Theta$.    
Similarly, we also have, $\Lambda^m(P^m(z)g(z))=0$ 
iff $f^m(\xi, z)g(z) \in \im\Theta$. From these two results, 
it is easy to see that the statements $(a)$ and $(b)$ 
in the theorem are equivalent to each other.
\epfv

Next we give two relations between the \JC and the {\bf IC}. 
The first one is as follows.

\begin{theo}\label{JC2IC}
The following statements are equivalent to each other.

$(a)$ The \JC holds for all $n\ge 1$.

$(b)$ For any $n\ge 1$ and homogeneous $P(z)\in \cz$ of degree $4$, the \IC holds for 
$f(\xi, z)=(\sum_{i=1}^n \xi_i^2)P(z)$ and $g(z)=P(z)$.
\end{theo} 

\pf By Theorem \ref{VC2JC}, we know that the \JC holds for all 
$n\ge 1$ iff for any $n\ge 1$ the \VC holds for the Laplace operator $\Delta=\sum_{i=1}^n\p_i^2$ and any homogeneous $P(z)=g(z)\in \bC[z]$ of degree $4$. Therefore, the equivalence of $(a)$ and $(b)$ in the theorem follows directly from this fact and 
Theorem \ref{VC2IC} above.
\epfv

The second relation between the \JC and the \IC is given by the following theorem.

\begin{theo}\label{AJC2IC}
The following statements are equivalent to each other.

$(a)$ The \JC holds for all $n\ge 1$.

$(b)$ For any $n\ge 1$ and homogeneous 
$H_i(z)\in \bC[z]$ $(1\le i\le n)$  of degree $3$, 
the \IC holds for $f(\xi, z)=\sum_{i=1}^n \xi_i H_i(z)$ 
and $g(z)=z_i$ for each $1\le i\le n$.
\end{theo} 

\pf First, by the classical homogeneous reduction 
(See \cite{BCW} and \cite{Y}) on the {\bf JC}, 
to prove or disprove the \JC holds, it is enough to consider polynomial maps $F(z)$ of the form $F(z)=z-H(z)$ with 
$H(z)\in \cz^{\times n}$ homogeneous 
of degree $3$.  

We fix a polynomial map $F(z)$ as above and 
denote by $G(z)$ the formal inverse map of $F(z)$.
Then by the Abhyankar-Gurjar inversion formula \cite{Ab} 
(See also \cite{BCW} and \cite{AG-VC}), we have
\begin{align}\label{AJC2IC-pe1}
 \sum_{m\ge 0} \sum_{\substack{\alpha \in \bN^n \\ |\alpha|=m}} \frac 1{\alpha!} \p^\alpha \left( H^\alpha(z) j(F)(z) 
g(z) \right)=g(G(z))
\end{align}
for any formal power series $g(z)\in \bC[[z]]$, where $j(F)$ denotes the {\it Jacobian} of the polynomial map $F$, i.e.\@ the determinant of the Jacobian matrix of $F$.

In particular, choose $g(z)=j(F)^{-1}\in \bC[[z]]$.
Since $j(F)(G)j(G) \equiv 1$, 
we get 
\begin{align}\label{AJC2IC-pe2}
\sum_{m\ge 0} \sum_{\substack{\alpha \in \bN^n \\ |\alpha|=m}} \frac 1{\alpha!} \p^\alpha \left( H^\alpha(z) \right)=j(F)^{-1}(G(z))=j(G)(z).
\end{align}

Since additionally $j(G)(F)j(F) \equiv 1$, we have that $j(F)(z)\equiv 1$ iff $j(G)(z)\equiv 1$. 
So by the equation above, we have 
\begin{align}\label{Equiv-0}
j(F)(z)\equiv 1 \Leftrightarrow 
\sum_{m\ge 0} \sum_{\substack{\alpha \in \bN^n \\ |\alpha|=m}} \frac 1{\alpha!} \p^\alpha \left( H^\alpha(z) \right)\equiv 1.
\end{align}

Since $\deg \p^\alpha \left( H^\alpha(z) \right)=2|\alpha|$ for any 
$\alpha\in \bN^n$, the equation above is equivalent to 
\begin{align}\label{AJC2IC-pe3}
\sum_{\substack{\alpha \in \bN^n \\ |\alpha|=m}} \frac 1{\alpha!} \p^\alpha \left( H^\alpha(z) \right)=0
\end{align}
for any $m\ge 1$.

Let $f(\xi, z):=\sum_{i=1}^n \xi_i H_i(z)$. 
Then for any $m\ge 1$, we have 
\allowdisplaybreaks{
\begin{align}
\mcE(f^m(\xi, z))=\mcE \left( \Big(\sum_{i=1}^n \xi_i H_i(z)\Big)^m\right) 
=\sum_{\substack{\alpha \in \bN^n \\ |\alpha|=m}} \frac {m!}{\alpha!} 
\mcE (\xi^\alpha H^\alpha(z)). \nno 
\end{align}
Applying Eq.\,(\ref{NewDef-E}) gives:
\begin{align}
\mcE(f^m(\xi, z))
=m! \sum_{\substack{\alpha \in \bN^n \\ |\alpha|=m}} \frac 1{\alpha!} \p^\alpha  (H^\alpha(z)). \label{AJC2IC-pe4} 
\end{align} }

By Eqs.\,(\ref{Equiv-0})--(\ref{AJC2IC-pe4}) and also 
Theorem \ref{Main-Thm-1}, we see that 
\begin{align}\label{Equiv-1}
j(F)(z)\equiv 1 & \Leftrightarrow f^m(\xi, z) \in \im\Theta 
\mbox{ for each $m\ge 1$}.
\end{align}

Now assume $j(F)\equiv 1$. By applying Eq.\,(\ref{AJC2IC-pe1}) 
to $g(z)=z$, we get 
\begin{align}\label{AJC2IC-pe4.2}
\sum_{m\ge 0} \sum_{\substack{\alpha \in \bN^n \\ |\alpha|=m}} \frac 1{\alpha!} \p^\alpha \left( H^\alpha(z)\, z \right)=G(z).
\end{align}

Note that, for any $m\ge 1$ and any $1\le i\le n$, we have  
\begin{align*}
\deg \sum_{\substack{\alpha \in \bN^n \\ |\alpha|=m}} \frac 1{\alpha!} \p^\alpha \left( H^\alpha(z) z_i \right)=2m+1.
\end{align*}

Therefore, from Eq.\,(\ref{AJC2IC-pe4.2}), we see that 
\begin{align}\label{AJC2IC-pe5}
\mbox{the \JC holds for $F(z)$} \Leftrightarrow 
 \sum_{\substack{\alpha \in \bN^n \\ |\alpha|=m}} \frac 1{\alpha!} \p^\alpha \left( H^\alpha(z) \, z \right)=0 \mbox{ when $m\gg 0$.}
\end{align}

But, on the other hand, by a similar argument as for 
Eq.\,(\ref{AJC2IC-pe4}), for any $m\ge 1$ and $1\le i\le n$, 
we also have
\begin{align}\label{AJC2IC-pe6}
\mcE(f^m(\xi, z)z_i) &=m! 
\sum_{\substack{\alpha \in \bN^n \\ |\alpha|=m}} 
\frac 1{\alpha!} \p^\alpha  (H^\alpha(z)\, z_i). 
\end{align}

Therefore, by  Eqs.\,(\ref{AJC2IC-pe5}), (\ref{AJC2IC-pe6}) and  Theorem \ref{Main-Thm-1}, we have that 
\begin{align}\label{Equiv-2}
\mbox{\it the }&\mbox{\it \JC holds for $F(z)$ } \\
& \Leftrightarrow f^m(\xi, z)z_i \in \im\Theta 
\mbox{ \it for each $1\le i\le n$ 
when $m\gg 0$.} \nno
\end{align}

Finally, from Eqs.\,(\ref{Equiv-1}) and (\ref{Equiv-2}), it is easy to see that the statements $(a)$ and $(b)$ in the theorem are indeed equivalent to each other.
\epfv

From the proofs of Theorems \ref{JC2IC} and \ref{AJC2IC}, we have the following ``non-trivial" families of polynomials all whose powers lie in $\im\Theta$.

\begin{corol}\label{Examples}
$(a)$ For any Hessian nilpotent polynomial 
$P(z)\in \cz$, i.e.\@  the Hessian matrix 
$\mbox{\rm Hes\,}(P)=\left(\frac{\p^2 P}{\p z_i \p z_j}\right)$ is nilpotent,
set $f(\xi, z)\!:=(\sum_{i=1}^n \xi_i^2)P(z)$. Then  
for any $m\ge 1$, we have $f^m(\xi, z)\in \im\Theta$. 

$(b)$ For any homogeneous 
$H=(H_1, H_2, ..., H_n)\in \cz^{\times n}$ of degree $d\ge 2$ with the Jacobian matrix $JH$ nilpotent, 
set $f(\xi, z)\!:=\sum_{i=1}^n \xi_i H_i(z)$. 
Then for any $m\ge 1$, we have $f^m(\xi, z)\in \im\Theta$. 
\end{corol}
\pf $(a)$ First, from the proof of Theorem \ref{JC2IC}, 
we see that for any $m\ge 1$, 
$f^m(\xi, z)\in \im\Theta$ iff $\Delta^m(P^m(z))=0$.
Second, by Theorem $4.3$ in \cite{HNP}, 
we know that $P(z)$ is Hessian nilpotent 
iff $\Delta^m(P^m(z))=0$ for all $m\ge 1$. 
Hence $(a)$ follows.

$(b)$ From the proof of Theorem \ref{AJC2IC}, we see that the statement holds if $H$ is homogeneous of degree $3$. 
But the argument there does not depend on the specific degree 
of $H$ but the homogeneity of $H$, as long as $d\ge 2$. Therefore, $(b)$ also follows from the same argument. 
\epfv

\subsection{Another Description of $\im\Theta$}
\label{S3.3}

In this subsection, we give a new description (see Theorem 
\ref{2nd-description}) for the image of 
the differential operators 
$\Theta=(\Theta_1, \Theta_2, ..., \Theta_n)$ 
with $\Theta_i=\xi_i-\p_i$ as before. The advantage of this description is that it does not involve any differential 
operators. We will give two different proofs for this result. 
While the first one is more straightforward and shorter, the second one is more conceptual. From the second proof, we will also see that 
the \IC is actually connected with integrals of polynomials 
such as multidimensional Laplace transformations of polynomials. 
For more discussions on the connections of the \IC with 
integrals of polynomials over open subsets of $\bR^n$, 
see \cite{GIC}.

First, let us fix the following notation and terminology.

For any Laurent polynomial $q(z)\in \czz$ and $\alpha\in \bZ^n$, we denote by $[z^\alpha]q(z)$ the coefficient of $z^\alpha$ of $q(z)$. For any subspace 
$V\subset \bC[z^{-1}, z]$, we say the $V$-part of $q(z)$ is zero if, for any 
$\alpha\in \bZ^n$ with $z^\alpha\in V$, we have 
$[z^\alpha]q(z)=0$. In particular, we also say that the holomorphic part of $q(z)$ 
is zero if the $\cz$-part of $q(z)$ is zero. 

We define a linear map $\mcZ: \cxz\to \czz$ 
by setting 
\begin{align}\label{Def-mcZ}
\mcZ(\, g(\xi) z^\beta\, ):=\beta!\, g(z^{-1})\,z^\beta
\end{align}
for each $g(\xi)\in \bC[\xi]$ and $\beta\in\bN^n$.

The main result of this subsection is the following theorem.

\begin{theo}\label{2nd-description}
For any $f(\xi, z)\in \cxz$,   
 $f(\xi,z)\in \im\Theta$ 
iff the holomorphic part of the Laurent polynomial  
$\mcZ(f(\xi, z))$ is equal to zero.
\end{theo}

\underline{\it First Proof:} 
We write  
$f(\xi, z)=\sum_{\alpha, \beta\in \bN^n} c_{\alpha, \beta} 
\xi^\alpha z^\beta$ for some $c_{\alpha, \beta}\in \bC$.
By Eq.\,(\ref{NewDef-E}), we have  
\begin{align*}
\mcE(f(\xi, z))&=\sum_{\alpha, \beta\in \bN^n} 
c_{\alpha, \beta} \p^\alpha z^\beta 
=\sum_{\substack{\alpha, \beta\in \bN^n\\ \alpha\le \beta}} 
\frac{\beta!}{(\beta-\alpha)!}\, c_{\alpha, \beta} 
 z^{\beta-\alpha} \\
&=\sum_{\gamma\in \bN^n}
\sum_{\substack{\alpha, \beta \in \bN^n \\ 
\beta-\alpha=\gamma} } \frac{\beta!}{\gamma!} \, 
c_{\alpha, \beta} z^{\gamma}
\end{align*}

On the other hand, by Eq.\,(\ref{Def-mcZ}), we have
\begin{align*}
\mcZ(f(\xi, z))=\sum_{\alpha, \beta\in \bN^n}
\beta! c_{\alpha, \beta} 
z^{\beta-\alpha}=\sum_{\gamma\in \bZ^n}
\sum_{\substack{\alpha, \beta \in \bN^n \\ 
\beta-\alpha=\gamma}} \beta! c_{\alpha, \beta} 
z^{\gamma}.
\end{align*}

By the two equations above, we see that 
for any $\gamma\in \bN^n$, we have 
\begin{align}
\gamma!\,[z^\gamma]\mcE(f(\xi, z))=[z^\gamma]\mcZ(f(\xi, z)).
\end{align}
Therefore, we have that $f(\xi, z)\in \Ker\mcE$ iff 
the holomorphic part of the Laurent polynomials  
$\mcZ(f(\xi, z))$ is equal to zero. By Theorem 
\ref{Main-Thm-1}, $\im\Theta=\Ker\mcE$, 
hence the theorem follows. 
\epfv

Next we give another proof for Theorem \ref{2nd-description} by using the {\it multi-dimensional Laplace transformations}. In order to do that, we first need to fix the following notation.

Let $\bR^+$ denote 
the set of all positive real numbers.  
By the restriction, we can and will view $f(\xi, z)$ as a family of $\bC$-valued functions in $z\in (\bR^+)^n$ parameterized by 
$\xi \in  (\bR^+)^n$, which we will still denote 
by $f(\xi, z)$. 

Now we consider the following integral:
\begin{align}\label{LaplaceTrans}
\cL(f)(\xi)\!:=\int_{(\bR^+)^n} f(\xi, z)e^{-\xi z}\, dz,
\end{align}
where $\xi z=\sum_{i=1}^n \xi_iz_i$.

Note that when $f(\xi, z)=h(z)$ 
for some $h(z)\in \cz$, i.e.\@  $f(\xi, z)$ is independent 
on $\xi$, $\cL(f)$ is nothing but the multidimensional Laplace 
transformation of the polynomial $h(z)$. In general, the integral in Eq.\,(\ref{LaplaceTrans}) is the base-field extension of the multidimensional Laplace transformation from $\cz$ to 
$\bC[\xi][z]=\cxz$. So we will still refer 
$\cL(f)$ as the ({\it multidimensional}) 
{\it Laplace transformation} of 
$f(\xi, z)\in \cxz$. \\

\underline{\it Second Proof of Theorem \ref{2nd-description}:}
We start with the following two observations on 
the Laplace transformation  
defined in Eq.\,(\ref{LaplaceTrans}).

First, for any $\beta\in \bN^n$, we have
\begin{align}\label{2nd-pf-pe1}
\cL(z^\beta)(\xi)&=
\int_{(\bR^+)^n} z^\beta e^{-\xi z}\, d z  
=\int_{(\bR^+)^n} (-\p_\xi)^\beta  e^{-\xi z}\, d z \\ 
&=(-\p_\xi)^\beta \int_{(\bR^+)^n}   e^{-\xi z}\, d z  
=(-\p_\xi)^\beta (\xi^{[-1]}) \nno\\
&=\beta! \, \xi^{-\beta} \xi^{[-1]}, \nno
\end{align}  
where $\xi^{[-1]}:=\prod_{i=1}^n \xi_i^{-1}$.

Therefore, for any polynomial 
$h(z)=\sum_{\alpha\in \bN^n} c_\alpha z^\alpha$ 
with $c_\alpha\in \bC$, we have
\begin{align} \label{2nd-pf-pe2}
\cL(h)(\xi) 
=\xi^{[-1]} \sum_{\alpha\in \bN^n} \alpha! c_\alpha \xi^{-\alpha} 
\in \xi^{[-1]}\bC[\xi^{-1}]. 
\end{align}  

Also, from the equation above, we see that  
$h(z)\neq 0$ iff $\cL(h)(\xi)\neq 0$.

Second, for any $1\le i\le n$ and 
$h(\xi, z)\in \cxz$, we have 
\begin{align*}
\int_0^\infty  e^{-\xi z} (\Theta_i h(\xi, z)) \, d z_i&=
-\int_0^\infty   e^{-\xi z} 
\Big( (\p_i-\xi_i) h(\xi, z)\Big) \, d z_i \\
&=-\int_0^\infty 
\p_i\left( h(\xi, z) e^{-\xi z}\right) \, d z_i \\
&=-\left. (h(\xi, z) e^{-\xi z}) \right|_{z_i=0}^{+\infty} 
=\left. ( h(\xi, z)e^{-\xi z} )\right|_{z_i=0}.
\end{align*} 

Hence we also have
\begin{align*}
\cL(\Theta_i h)(\xi)=\int_{(\bR^+)^{n-1}}
\left. ( h(\xi, z)e^{-\xi z} )\right|_{z_i=0} 
dz_1\cdots dz_{i-1}dz_{i+1}\cdots dz_n.
\end{align*}

From the equation above, it is easy to see that 
$\cL(\Theta_i h)(\xi)$ can not have any $\xi^\alpha$-term with the $i^{\rm th}$ component of $\alpha$ strictly less than $0$.  Therefore, for any $g(\xi, z)\in \im\Theta$, its Laplace transformation $\cL(g)(\xi)$ can not have any non-zero 
$\xi^{[-1]}\bC[\xi^{-1}]$ part.

Now, for any $f(\xi, z)\in \cxz$, by Theorem \ref{D-Taylor} 
we may write $f(\xi, z)$ uniquely as 
$f(\xi, z)=a(z)+g(\xi, z)$ with $a(z)\in \cz$ 
and $g(\xi, z)\in \im\Theta$. Then 
by the observations above, we see that
the $\xi^{[-1]}\bC[\xi^{-1}]$-part of 
$\cL(f)(\xi)$ is equal to $\cL(a)(\xi)$.  
Therefore, we have that  
$f(\xi, z)\in \im\Theta$ iff $a(z)=0$ 
(by Theorem  \ref{D-Taylor})
iff $\cL(a)=0$ iff the $\xi^{[-1]}\bC[\xi^{-1}]$-part of 
$\cL(f)(\xi)$ is equal to zero.

On the other hand, if we write $f(\xi, z)=\sum_{\beta\in \bN^n} b_\beta(\xi)z^\beta$ with $b_\beta(\xi)\in \bC[\xi]$. Then by 
Eqs.\,(\ref{2nd-pf-pe1}) and (\ref{Def-mcZ}), we have 
\begin{align}\label{2nd-pf-pe3}
\cL(f)(\xi) &=\sum_{\beta\in \bN^n}
b_\beta(\xi) \int_{(\bR^+)^n} z^\beta e^{-\xi z}\, d z \\
&=\xi^{[-1]}\sum_{\beta\in \bN^n}  \beta! \, b_\beta(\xi) \xi^{-\beta} 
=\xi^{[-1]}\mcZ(f)(z)|_{z=\xi^{-1}}.\nno
\end{align}  
Therefore, we have that the $\xi^{[-1]}\bC[\xi^{-1}]$-part of $\cL(f)(\xi)$ is equal to zero iff the $\bC[\xi^{-1}]$-part of 
$\mcZ(f)(z)|_{z=\xi^{-1}}$ is equal to zero iff  
the holomorphic part of 
$\mcZ(f)(z) \in \czz$ is equal to zero. 
Hence the theorem follows. 
\epfv

From the second proof above, it is easy to see that Theorem 
\ref{2nd-description} can be re-stated in terms of multi-dimensional Laplace transformations as follows.

\begin{corol}\label{L2nd-description}
For any $f(\xi, z)\in \cxz$,   
 $f(\xi,z)\in \im\Theta$ iff the $\xi^{[-1]}\bC[\xi^{-1}]$-part of its Laplace transformation $\cL(f)(\xi)$ is 
equal to zero.
\end{corol}

By Theorem \ref{2nd-description} and Corollary 
\ref{L2nd-description} above, we see that 
the \IC can be re-stated as follows.

\begin{conj}\label{RestateIC}
Let $\cM$ be the subspace of polynomials 
$f(\xi, z)\in \cxz$ such that $\mcZ(f)$ has no 
holomorphic part, or equivalently, the Laplace 
transformation $\cL(f)(\xi)$ has no 
$\xi^{[-1]}\bC[\xi^{-1}]$-part. Then for any 
$f(z), g(z)\in \cxz$ with $f^m\in \cM$ for 
each $m\ge 1$, we have $f^m g\in \cM$ 
when $m\gg 0$.
\end{conj}

For comparison, let us point out that
the following theorem first conjectured 
in \cite{GVC} has recently been proved in 
\cite{EWZ} by using 
the Duistermaat-van der Kallen theorem, Theorem \ref{ThmDK}.

\begin{theo}
Let $\cM$ be the subspace of Laurent polynomials 
$f(z)\in \czz$ such that $f(z)$ has no holomorphic part.  
Then for any $f, g\in \czz$ with $f^m\in \cM$ for 
each $m\ge 1$, we have $f^m g\in \cM$ when $m\gg 0$.
\end{theo}

Actually, as pointed out in \cite{GVC} and \cite{EWZ}, 
the theorem above implies that the \VC holds 
when either $\Lambda$ is a monomial of 
$\p$ or $P(z)$ is a monomial of $z$ 
(See \cite{GVC} and \cite{EWZ} for more details). 
Consequently, by Theorem \ref{VC2IC}, we see that  
the \IC holds when $f(z, \xi)\in \cxz$ has 
the form $\xi^{\alpha}P(z)$ or $\Lambda(\xi)z^\alpha$ 
with $\Lambda(\xi)\in \cxi$, $P(z)\in \cz$ 
and $\alpha\in \bN^n$.

{\small \sc Department of Mathematics, Illinois State University,
Normal, IL 61790-4520.}

{\em E-mail}: wzhao@ilstu.edu.

\end{document}